\documentclass[10pt]{amsart}
\usepackage{amscd}
\usepackage{xypic}  
\usepackage{amssymb}
\usepackage{amsthm}
\usepackage{epsfig}
\usepackage{paralist}
\usepackage{comment}
\usepackage[T1]{fontenc}
\usepackage[all,cmtip]{xy}
\usepackage{amsmath}
\usepackage{color}

\def\bdi{\begin{diagram}}
\def\edi{\end{diagram}}
\def\bsit{\begin{sit}}
\def\esit{\end{sit}}
\def\brem{\begin{rem}}
\def\brems{\begin{rems}}
\def\erem{\end{rem}}
\def\erems{\end{rems}}
\def\bprob{\begin{prob}}
\def\eprob{\end{prob}}
\def\bprobs{\begin{probs}}
\def\eprobs{\end{probs}}
\def\bques{\begin{ques}}
\def\eques{\end{ques}}
\def\bexa{\begin{exa}}
\def\bexas{\begin{exas}}
\def\eexa{\end{exa}}
\def\eexas{\end{exas}}
\def\bdefi{\begin{defi}}
\def\edefi{\end{defi}}
\def\bdefis{\begin{defis}}
\def\edefis{\end{defis}}
\def\bcor{\begin{cor}}
\def\ecor{\end{cor}}
\def\blem{\begin{lem}}
\def\elem{\end{lem}}
\def\bconv{\begin{conv}}
\def\econv{\end{conv}}
\def\bconj{\begin{conj}}
\def\econj{\end{conj}}
\def\bprop{\begin{prop}}
\def\eprop{\end{prop}}
\def\bthm{\begin{thm}}
\def\ethm{\end{thm}}
\def\bnota{\begin{nota}}
\def\enota{\end{nota}}
\def\bsit{\begin{sit}}
\def\esit{\end{sit}}
\def\be{\begin{equation}}
\def\ee{\end{equation}}
\def\bproof{\begin{proof}}
\def\eproof{\end{proof}}
\def\ba{\begin{array}}
\def\ea{\end{array}}
\def\la{\label}
\def\bcl{\begin{claim}}
\def\ecl{\end{claim}}

\newtheorem{theorem}{Theorem}[section]

\theoremstyle{definition}

\theoremstyle{remark}

\numberwithin{equation}{section}

\newtheorem{cor}[theorem]{Corollary}
\newtheorem{lem}[theorem]{Lemma}
\newtheorem{prop}[theorem]{Proposition}
\newtheorem{claim}[theorem]{Claim}
\newtheorem{defi}[theorem]{Definition}
\newtheorem{defis}[theorem]{Definitions}
\newtheorem{conj}[theorem]{Conjecture}

\newtheorem{conv}[theorem]{Convention}
\newtheorem{nota}[theorem]{Notation}
\newtheorem{rem}[theorem]{Remark}
\newtheorem{rems}[theorem]{Remarks}
\newtheorem{exa}[theorem]{Example}
\newtheorem{exas}[theorem]{Examples}

\newcommand{\Z}{\ensuremath{\mathbb{Z}}}

\newcommand{\cO}{{\ensuremath{\mathcal{O}}}}

\def\CC{{\mathbb C}}

\def\PP{{\mathbb P}}
\def\Of{{\mathcal{O}}}

\def\deg{\mathop{\rm deg}}

\def\Pic{\mathop{\rm Pic}}

\renewcommand{\epsilon}{\varepsilon}
\renewcommand{\phi}{\varphi}

\newcommand{\bnum}{\begin{enumerate}}
\newcommand{\enum}{\end{enumerate}}

\begin{document}
\title[A remark on the intersection of plane curves]{A remark on the intersection of plane curves}

\author{C.~Ciliberto}
\address{Dipartimento di Matematica, Universit\`a degli
Studi di Roma ``Tor Vergata'', Via della Ricerca Scientifica,
00133 Roma, Italy}
\email{cilibert@mat.uniroma2.it}

\author{F.~Flamini}
\address{Dipartimento di Matematica, Universit\`a degli
Studi di Roma ``Tor Vergata'', Via della Ricerca Scientifica,
00133 Roma, Italy}
\email{flamini@mat.uniroma2.it}

\author{M.~Zaidenberg}
\address{Universit\'e Grenoble Alpes, Institut Fourier,
CS 40700, 38058 Grenoble cedex 09, France}
\email{Mikhail.Zaidenberg@univ-grenoble-alpes.fr}


\thanks{{\bf Acknowledgements:} This research started during a visit of the third author at the Dept. of Mathematics of Univ. Roma ``Tor Vergata'' (supported by this Department, the INdAM ``F. Severi" in Rome and the cooperation program GDRE-GRIFGA) and of the second author at the Institut Fourier (supported by INdAM-GNSAGA and within the context of the Laboratory Ypatia of Mathematical Sciences LIA LYSM AMU-CNRS-ECM-INdAM). The authors  thank all these institutions and programs for the support and the excellent working conditions. The first author acknowledges the {\em MIUR Excellence Department Project} awarded
to the Department of Mathematics, University of Rome Tor Vergata, CUP E83C18000100006. The work of the third author was partially supported by the grant
346300 for {\em IMPAN} from the {\em Simons Foundation} and the matching 2015-2019 Polish
MNiSW fund (code: BCSim-2018-s09). Collaboration has also benefitted of funds {\em Mission Sustainability 2017 - Fam Curves} CUP E81-18000100005 (Tor Vergata University)}

\begin{abstract} Let $D$ be a very general curve of  degree $d=2\ell-\epsilon$ in $\PP^2$, with $\epsilon\in \{0,1\}$. Let  $\Gamma \subset \PP^2$ be an integral curve of geometric genus $g$  and degree $m$, $\Gamma \neq D$, and let $\nu: C\to \Gamma$ be the normalization. Let  $\delta$ be the degree of the \emph{reduction modulo 2} of the divisor $\nu^*(D)$ of $C$ (see \S\;\ref{ss:delta}).  In this paper we prove the inequality 
$4g+\delta\geqslant m(d-8+2\epsilon)+5$. We compare this with similar inequalities due to Geng Xu (\cite{Xu1,Xu2}) and Xi Chen (\cite{chen1,chen2}).  Besides, we provide a brief account on genera of subvarieties in projective hypersurfaces.
\end{abstract}

\maketitle

\tableofcontents



\section*{Introduction} Given an effective divisor $D\in |\Of_{\PP^n}(d)|$  and an integral (i.e., reduced and
irreducible) projective curve $\Gamma$ of degree $m$ in $\PP^ n$, which is not contained in supp$(D)$,
let 
$j(D,\Gamma)$ be the order of $\Gamma\cap D$. Assume $D$ is very general and set
$$j(n,d,m):= \min\{j(D,\Gamma)\,|\, \Gamma \subset \PP^n  \; \mbox{as above}\}\quad \text{and}\quad j(n,d):= \min_{m\geqslant 1}\{j (n,d,m)\}.$$ Similarly, with $\Gamma$ and $D$ as before, let $i(D,\Gamma)$ stand for the {\em number of places} of $\Gamma$ on $D$, that is, the number of centers of local branches of the curve $\Gamma$ on $D$. Then, set
$$i(n,d,m):= \min\{i(D,\Gamma)\,|\, \Gamma \subset \PP^n \; \mbox{as above}\}\quad \text{and}\quad i(n,d):= \min_{m\geqslant 1}\{i (n,d,m)\}.$$ 

The problem of computing $j(n,d)$ and $i(n,d)$ has been considered in \cite{chen1, Xu1,Xu2} (basically devoted to $n=2$ case) and \cite{chen2} (where the case $n\geqslant 2$ is considered).  The relations of this with  the famous Kobayashi problem on hyperbolicity  of  the complement of a very general hypersurface in $\PP^ n$ is well known and we do not dwell on this here (see, e.g., 
\cite {chen2}).

Geng Xu (\cite[Thm.\,1]{Xu1}) proved that $$j(2,d)=d-2,\,\;\mbox{for any}\; d\geqslant 3,$$where the equality is attained either by a bitangent line or by an inflectional tangent line of $D$, i.e. the minimum is achieved by $m=1$. Moreover, for $d=3$, he also proved in \cite[Corollary]{Xu2} that, for any integer $m \geqslant 1$, the number of rational curves of degree $m$ which meet set-theoretically a given (arbitrary) smooth plane cubic curve $D$ at
exactly one point is finite and positive. Therefore,
for $d=3$ the minimum $j(2,3)=1$ is achieved by any integer $m \geqslant 1$.  

Xi Chen (\cite[Thm.\,1.2]{chen1}) proved  that, for $d>m$, one has
$$j(2, d, m) \geqslant \min \left\{dm - \frac{m(m + 3)}{2}, \; 2dm - 2m^2 - 2\right\}.$$ {Furthermore (cf. \cite[Cor.\,1.1]{chen1}),  for $d \geqslant {\rm max} \{\frac{3m}{2} - 1, \; 3\}$ one has
$$j(2,d, m) = dm -   \dim (|\mathcal{O}_D(m)|)=dm- \frac{m(m + 3)}{2}\,\, .$$
In addition, he conjectured (see \cite[Conj.\,1.1]{chen1}) 
that 
$$j(2,d, m) = dm -   \dim (|\mathcal{O}_D(m)|) \quad \text{if}\quad d>\max\{m,2\}.$$

In arbitrary dimension $n\geqslant 2$, Xi Chen (\cite[Thm\;1.7]{chen2})  proved that, for  $D$ very general and $\Gamma$ as above, one has
\be\la{eq:Chen} 2g -2+ i(D,\Gamma)\geqslant (d-2n) m \,,\ee
where $g$ is the \emph{geometric genus} of $\Gamma$, i.e., the arithmetic genus of its normalization.

In this paper we obtain a new inequality of  type \eqref{eq:Chen}, although only in the case $n=2$ (see Theorem \ref {prop:even}).  Indeed, let $D$ be a very general curve of  degree $d=2\ell-\epsilon$ in $\PP^2$, with $\epsilon\in \{0,1\}$. Let  $\Gamma$ be an integral curve in $\PP^2$ of geometric genus $g$  and degree $m$, $\Gamma \neq D$, and let $\nu: C\to \Gamma$ be the normalization. Let  $\delta(D,\Gamma)$ be the degree of the \emph{reduction modulo 2} of the divisor $\nu^*(D)$ on $C$ (cf.\;\S\;\ref{ss:delta}). 
In Theorem \ref {prop:even} we prove that 
\begin{equation}\label{eq:us}
4g+\delta(D,\Gamma)\geqslant m(d+2\epsilon-8)+5.
\end{equation}
Note that $\delta(D,\Gamma)\leqslant i(D,\Gamma)$, and  the equality  holds if and only if at any place  $p$ of 
$\Gamma$ on $D$, the local intersection multiplicity of $D$ and $\Gamma$ at $p$ is odd. This happens, for instance, if $\Gamma$ intersects $D$ transversely. In the latter case $\delta(D,\Gamma)= i(D,\Gamma)=md$ and both  \eqref {eq:Chen} and \eqref {eq:us} are uninteresting. On the other hand, \eqref {eq:Chen} and \eqref {eq:us} become interesting when $\delta(D,\Gamma)$ and $i(D,\Gamma)$ are small. Though the difference between the two quantities is a priori unpredictable, one may expect that, generally speaking, $\delta(D,\Gamma)$ is strictly smaller than $i(D,\Gamma)$. Unfortunately, the genus $g$ works against us in \eqref {eq:us}; however, for $g=0,1$ and $d$ even, \eqref {eq:us} is better than \eqref{eq:Chen}. 
Further related problems have been recently considered in \cite{ChenZhu,Liu1,Liu2}. 

As a final remark, note that \eqref {eq:us} is more useful than \eqref {eq:Chen} if one looks, as we do in this paper, at the geometric genera of curves contained in a \emph{double plane} $X_d$, that is, a cyclic double cover of $\PP^2$ branched along a very general plane curve $D$ of even degree $d$. For instance, letting $g=0$, $\delta(D,\Gamma)=0,2$ and $d$ even, we are looking actually for rational curves on $X_d$. By \eqref {eq:us} we see that such a rational curve over $\Gamma$ might exist, as expected, only for $d\leqslant 6$ (for low $m$ one has even smaller bounds on $d$).  The case $d=6$ corresponds to a $K3$ surface, which always contains infinitely many rational curves. In contrast, it follows from  \eqref {eq:us} that the double planes with very general branching curves of even degree $d\geqslant 8$ ($d\geqslant 10$, respectively) do not carry any rational curve (any rational or elliptic curve, respectively, hence are \emph{algebraically hyperbolic}). For $d=8$ and $d=10$ these double  planes are \emph{Horikawa surfaces} $H_8$ and $H_{10}$, that is, their Chern numbers satisfy $c_2 = 5c^2_1 + 36$ (in other words, $(c_1^2,c_2)$ lies on the Noether line). The algebraic hyperbolicity of $H_{10}$ was established first by X.~Roulleau and E.~Rousseau (\cite{RR}). J.\ Liu (\cite{Liu2}) showed that some of the Horikawa surfaces $H_{10}$ are even  Kobayashi hyperbolic, whereas there is  no hyperbolic $H_8$. Indeed, the Horikawa double planes $H_8$ carrying elliptic curves are dense in the set of all such surfaces, while the Kobayashi hyperbolicity is open in the Hausdorff topology.

The proof of Theorem \ref {prop:even} presented in \S \ref {sec:result} follows, with minor variations due to the different setting, the basic ideas exploited in \cite{CL} (and later in \cite {CLR}). These are based on a smart use of the theory of {\em focal loci}, see e.g. \cite{CC}. For the reader's convenience, we recall   in \S\,\ref{App} the basic notions and results of this theory.  We apply this technique to families of double covers of $\PP^2$ branched along a very general plane curve $D$ or along $D$ plus a general line, according to whether the degree of $D$ is even (see \S\;\ref{ss:even} and \S\;\ref{ss:proofeven}) or odd (see \S\;\ref {ss:odd} and \S\; \ref {ss:proofodd}).

In the last \S \ref{sec:appendix} we provide a short survey on genera of subvarieties in projective varieties, with accent on projective hypersurfaces. 


\section{Focal loci}\label{App} For the reader's convenience, we recall here some basic notions from \cite{CC, CL}.

Let $X$ be a smooth projective variety of dimension $n+1$. Assume we have a flat, projective family $\mathcal{D}\stackrel{p}{\longrightarrow}\mathcal{B}$  of effective divisors on $X$ over a smooth, irreducible, quasiprojective base
$\mathcal{B}$, with irreducible general fiber. Up to shrinking $\mathcal B$ to a suitable Zariski dense, open subset, we may suppose that for any closed point $b \in \mathcal B$ the fiber $D_b$ of $\mathcal{D}\stackrel{p}{\longrightarrow}\mathcal{B}$ over $b$ is irreducible. 

Assume we have a commutative diagram   
\begin{equation}\label{eq:in1}
\xymatrix{
\mathfrak C \ar[drr]_{\mathfrak q}  \quad \ar@{^{(}->}[rr] ^ i &&\,\,  \mathcal{D}\ar[d]^p \\
 &&  \,\, \mathcal B
 }
\end{equation}
where $\mathfrak q: \mathfrak {C}\to \mathcal B$ is a flat projective family such that, for all $b\in \mathcal B$, the fiber $\Gamma_b$ over $b$ is a reduced curve of \emph{geometric genus} $g$, and where  $i$ is an inclusion: so, for any $b \in \mathcal B$,  one has $\Gamma_b\subset D_b$ via the inclusion $i_b$. 

By a result of Tessier (see \cite[Th\'eor\`eme 1]{Teis}), there is a simultaneous normalization 
\be\label{eq:tess}
\xymatrix{
{\mathcal{C}}  \ar[drr]_{q}  \quad  \ar[rr]^ \nu  &&\,\, \mathfrak C\ar[d]^{\mathfrak q} \\
 &&  \,\, \mathcal B
 }
\ee
such that $\mathcal C$ is smooth and,  for every $b\in \mathcal B$, the fiber $C_b$ of $q: {\mathcal{C}} \to \mathcal B$ is the normalization $\nu_b: C_b\to \Gamma_b$ of $\Gamma_b$. For any $b\in \mathcal B$, the curve $C_b$ is smooth of \emph{(arithmetic) genus} $g$. 

Composing with the inclusion $\mathcal D \stackrel{j}{\hookrightarrow} \mathcal B \times X$, we get the commutative diagram 
\be\la{eq:bidg}
\xymatrix{
{\mathcal{C}}  \ar[drr]_{q}  \quad \quad  \ar[rr]^ \nu  &&\,\,\, \mathfrak C\ar[d]_{\mathfrak q} \quad \ar@{^{(}->}[rr] ^ i &&\,\,  \mathcal{D}\ar[d]^p\quad \ar@{^{(}->}[rr] ^ j&&\,\,  \mathcal B \times X \ar[dll]^{\text{pr}_1}  \ar[d]^{\text{pr}_2}\\
 &&  \,\, \mathcal B  \ar[rr]_{\text{id}}&& \mathcal B &&  X &&
 }
\ee
where $\text{pr}_i$ is the projection onto the $i$th factor, for $i=1,2$. 

We set $$s:=j\circ i \circ \nu \colon {\mathcal{C}}\to\mathcal{B}\times X\,,$$
and let  $\mathcal{N}:=\mathcal{N}_s$ be the normal sheaf to $s$, defined by the exact sequence 
\[ 
0 {\longrightarrow} \mathcal T_{\mathcal C} \stackrel{ds}{\longrightarrow}   s^*(\mathcal T_{\mathcal{B}\times X})  {\longrightarrow} \mathcal N{\longrightarrow} 0\,,
\]
where 
$\mathcal{T}_Y$ stands for the tangent sheaf of a smooth variety $Y$.  

For $b \in \mathcal B$ we set 
$$N_b:=\mathcal{N}|_{C_b}=\mathcal{N}\otimes\mathcal{O}_{C_b}\quad\mbox{and}\quad s_b=s|_{C_b}\colon C_b\to\{b\}\times X=X\,.$$

In addition, we let 
$$\phi:={\rm pr}_2 \circ s: \mathcal C\to X\,.$$
Then $\phi_b=\phi|_{C_b}$ coincides with $s_b$ for any $b\in \mathcal B$, that is,
\[ \xymatrix{ \phi_b=s_b:  \,\, C_b  \,\,  \ar[rr]^ {\nu_b}  &&\Gamma_b   \,\,  \ar@{^{(}->}[rr]^ {i_b} && D_b   \,\, \ar@{^{(}->}[rr] ^ {(\text{pr}_2\circ j)_b} && X }.\]
As in \cite[\S\,2]{CL}, we set \begin{equation}\label{eq:zetaC}
z({\mathcal{C}}):=  \dim\; (\phi(\mathcal{C})), 
\end{equation}so that $z({\mathcal{C}})\leqslant n+1=\dim (X).$ If $z({\mathcal{C}})=n+1$ one says that 
$\mathfrak C \stackrel{\mathfrak q}{\longrightarrow} \mathcal B$, or $\mathcal C\stackrel{q}{\longrightarrow} \mathcal B$, is a \emph{covering family}.

\bprop [See {\cite[Prop.\;1.4 and p. 98]{CC}}]\la{prop:CL} In the above setting, we have: 
 \begin{enumerate}
\renewcommand\labelenumi{\rm (\alph{enumi})}
\renewcommand\theenumi{\rm (\alph{enumi})}
\item
for any $b\in \mathcal B$, the sheaf $N_b$ fits into  the exact sequence
$$0\longrightarrow \mathcal {T}_{C_b}\stackrel{ds_b}{\longrightarrow} s_b^*(\mathcal{T}_X)\longrightarrow N_b\longrightarrow 0\,$$
and ${\mathcal{C}}\stackrel{q}{\longrightarrow}\mathcal{B}$ induces on $C_b$ a \emph{characteristic map} 
$$\chi_b\colon T_{ \mathcal B,b} \otimes\mathcal{O}_{C_b}\longrightarrow N_b\,,$$where $ T_{ \mathcal B,b} $ denotes the tangent space to $\mathcal B$ at $b$;
\item
 if $b\in \mathcal B$  and 
$x\in C_b$ are general points, 
then 
$$\dim \,(N_{b,x})=\dim \,(s_b^*(\mathcal{T}_X)_{x})-\dim \,(\mathcal{T}_{C_b, x})=n\quad  \text{and}\quad  {\rm rk}\, (\chi_{b,x})=z({\mathcal{C}})-1\,.$$
\end{enumerate}\eprop

\bdefi[See {\cite[Def.~ (1.5)]{CC}}]\la{def:focal-locus} Given $b\in \mathcal B$, the \emph{focal set} of  ${\mathcal{C}}\stackrel{q}{\longrightarrow}\mathcal{B}$ on $C_b$ is the closed subscheme $\Phi_b$ of $C_b$  defined as 
$$\Phi_b:=\{x\in C_b\vert \, {\rm rk}(\chi_{b,x})<z({\mathcal{C}})-1\}.$$
If $b\in \mathcal B$ is general, then $\Phi_b$ is a proper subscheme of $C_b$. 
The points in $\Phi_b$ are called \emph{focal points} of ${\mathcal{C}}\stackrel{q}{\longrightarrow}\mathcal{B}$ on $C_b$.
We  denote by $\Phi_b^{\rm sm}$ the open subset of $\Phi_b$ consisting of the points $x\in \Phi_b$ which map to smooth points of $\Gamma_b$ via the normalization morphism $\nu_b\colon C_b\to \Gamma_b$.
\edefi

\bprop
[{\cite[Prop.\;2.3 and Prop.\,2.4]{CL}}] \la{prop:CL-3} Let  ${\mathcal{C}}\longrightarrow\mathcal{B}$ be a covering family. Then the following hold.
\begin{enumerate}

\item
[{\rm (i)}] 
Suppose that for $x\in C_b$ the point
$s(x)$ is smooth in both $\Gamma_b$ and $D_b$. Assume also that $s(x)$ is a \emph{fundamental point} of the family $\mathcal{D}\stackrel{p}{\longrightarrow} \mathcal{B}$, i.e. it is a base point of the family. Then $x\in \Phi_b^{\rm sm}$.

\item 
[{\rm (ii)}] 
For a general point $b\in \mathcal B$ one has
\begin{equation}\label{eq:foc-set-deg} \deg (\Phi_b^{\rm sm}) \leqslant 2g-2-K_X\cdot \Gamma_b\,.\end{equation}
\end{enumerate} \eprop


\section{Double planes}\label{S:double} 
In this section we collect useful material for the proof of our main result. The result itself is stated and proven in \S \ref {sec:result}. The contents of this section, which suffice for our applications, can be easily adapted to the higher dimensional case.


\subsection{The $\delta$--invariant}\label{ss:delta} 
 Let $C$ be any smooth, irreducible, projective 
curve, and let $\Delta=\sum_i m_ip_i$ be an effective divisor on $C$. We set ${\Delta}_2:=\sum_i \overline m_ip_i\,$,  where $\overline m_i \in\{0,1\}$ is the residue of the integer $m_i$ modulo $2$. We also set $\delta_2(\Delta):=\deg (\Delta_2)$.

For any smooth curve $D \subset \PP^2$ and any integral curve $\Gamma \subset \PP^2$, $\Gamma\neq D$, with normalization $\nu\colon C\to \Gamma$, we set 
\begin{equation}\label{eq:delta}
\delta(D,\Gamma):=\delta_2(\nu^*(D))\,.
\end{equation} 
We notice that
$$
\delta(D,\Gamma)\leqslant i(D,\Gamma).
$$


\subsection{Basics on a certain weighted projective $3$-space}\label{ss:gengen} 
For any positive integer $\ell$, we denote by $\mathcal L_{\ell}$ the linear system $|\mathcal O_{\PP^ 2}(\ell)|$ of plane curves of degree $\ell$, and by $\mathcal U_{\ell}$ its open dense subset of points corresponding to smooth curves. We let $N_\ell=\dim(\mathcal L_{\ell})  =\frac {\ell(\ell+3)}2$. We denote by $\mathcal D_{\ell}\to \mathcal L_{\ell}$ the universal curve, and we  use the same notation $\mathcal D_{\ell}\to \mathcal U_{\ell}$ for its restriction to  $\mathcal U_{\ell}$. 

The linear system $\mathcal L_{\ell}$ determines the $\ell$th \emph{Veronese embedding} $\PP^2\stackrel {v_{\ell}} \hookrightarrow\PP^{N_{\ell}}$, whose image is the $\ell$--\emph{Veronese surface} $V_{\ell}$ in $\PP^{N_{\ell}}$. 
Let $[x]=[x_0,x_1, x_2]$ be homogeneous coordinates  in $\PP^ 2$, and let
$$[x^I], \,\,\, \text {where}\,\,\, I=(i_0,i_1,i_2) \,\,\, \text {is a multiindex such that} \,\,\, |I|=i_0+i_1+i_2=\ell, \,\,\, $$ 
be homogeneous coordinates
in $\PP^{N_{\ell}}$. In these coordinates
the {Veronese map} is given by
\[
\PP^2\ni [x] \stackrel{v_\ell}{\longrightarrow}[x^ I]_{|I|=\ell}\in \PP^{N_{\ell}},\,\,\, \text{where}\,\,\, x^ I:= x_0^ {i_0}x_1^ {i_1} x_2^ {i_2}.\]
We equip the weighted projective $3$-space $\PP(1,1,1,\ell)$ with  weigthed homogeneous  coordinates $[x,z]:=[x_0,x_1, x_2,z]$, where
$x_0,x_1, x_2$ [resp. $z$] have \emph{weigth} 1 [resp. has weight $\ell$]. We 
introduce as well coordinates  $[x^I,z]_{|I|=\ell}$ in $\PP^{N_{\ell}+1}$ and embed 
$\PP^{N_{\ell}}$ in $\PP^{N_{\ell}+1}$ as the hyperplane $\Pi$ with equation $z=0$. Then
 $\PP(1,1,1,\ell)$ can be identified with the cone $W_{\ell}\subset\PP^{N_{\ell}+1}$ over the $l$--Veronese surface $V_{\ell}$ with vertex $P=[0,\ldots,0,1]$.  Blowing $P$ up yields a minimal resolution 
 $$\rho: Z_{\ell}\to W_{\ell}\cong\PP(1,1,1,\ell),$$
 with exceptional divisor $E\cong V_{\ell}\cong \PP^2$. The projection from $P$  induces a $\PP^1$--bundle structure 
 $$\pi\colon Z_{\ell} \to V_{\ell}\cong \PP^2\,.$$Let $f$ be the class of a fiber of $\pi$. One has
$$Z_{\ell}\cong \PP(\mathcal O_{\PP^ 2}(\ell)\oplus \mathcal O_{\PP^ 2})\quad  \text{and}\quad \cO_{ Z_{\ell}}(1)=\rho^ *(\cO_{W_{\ell}}(1)).$$ 

For every integer $m$, we set
\begin{equation}\label{eq:lel}\mathcal O_{\ell}(m):=\pi^ *(\cO_{\PP^2}(m)) \quad \text{and}\quad \mathcal L_\ell(m):=|\mathcal O_{\ell}(m)|\,.\end{equation}
  Note that
\begin{equation}\label{eq:le}
\cO_{ Z_{\ell}}(1)\cong \mathcal O_{\ell}(\ell)\otimes \cO_{ Z_{\ell}}(E).
\end{equation}
Since 
$$E\cong\PP^2, \,\,\, \cO_{ Z_{\ell}}(1)\otimes \cO_E\cong \cO_E,\,\,\,  \text{and} \,\,\, \mathcal O_{\ell}(\ell)\otimes \cO_E\cong \cO_{\PP^ 2}(\ell),$$
we deduce
\begin{equation}\label{eq:kappa}
\cO_{ Z_{\ell}}(E) \otimes \mathcal O_E\cong \cO_{\PP^ 2}(-\ell).
\end{equation}

Finally, we denote by $K_{\ell}$ the canonical sheaf of $Z_{\ell}$. 
 
\blem\label{lem:canclass} One has 
$$K_{\ell}\cong \mathcal O_{\ell}(\ell-3)\otimes \cO_{ Z_{\ell}}(-2)\cong \mathcal O_{\ell}(-\ell-3)\otimes \cO_{ Z_{\ell}}(-2E)\,.$$
\elem

\bproof The Picard group $\Pic (Z_{\ell})$ is freely generated by the classes $\mathcal O_{\ell}(1)$ and $\cO_{ Z_{\ell}}(E)$, and also  by $\mathcal O_{\ell}(1)$ and $\cO_{ Z_{\ell}}(1)$, see \eqref {eq:le}.
Let $H$ [resp.~ $L$] be a general member of
$|\cO_{ Z_{\ell}}(1)|$ [resp.~ of 
$\mathcal L_{\ell}(1)$].  Write $K_{\ell}\sim\alpha H+\beta L$, where $\alpha,\beta\in\Z$. From the relations $K_{\ell}\cdot f=-2$, $H\cdot f=1$, and $L\cdot f=0$  one gets $\alpha=-2$. 

By  adjunction formula and \eqref {eq:kappa} we obtain 
$$\mathcal{O}_{\PP^2}(-3)\cong K_E=(K_{\ell}+E)|_E\cong (-2H+\beta L+E)|_E\cong\mathcal{O}_{\PP^2}(\beta-\ell)\,.$$ So, $\beta=\ell-3$, as desired. 
\eproof 

Finally, let $\mathbb G_{\ell}$ be the group of all  automorphisms of $\PP(1,1,1,\ell)$ which stabilize the divisor with equation $z=0$. This group is naturally isomorphic to the automorphism group of the pair $(W_{\ell},V_\ell)$, i.e. automorphisms of $W_{\ell}$ stabilizing $V_{\ell}$, where $V_\ell$ is cut out on $W_\ell$ by $\Pi$. In turn, the latter group is isomorphic to the automorphism group of the pair  $(Z_{\ell}, \rho^*(V_\ell))$. One has the exact sequence
\be\label{eq:gr}
1\to \mathbb C^ *\to \mathbb G_{\ell}\to {\rm PGL}(3,\mathbb C)\to 1.\ee


\subsection{The even degree case}\label{ss:even} 
Let $D$ be a smooth curve in $\PP^2$ of even degree $d=2\ell \geqslant 2$ which, in the homogeneous coordinate system fixed in \S\ref{ss:gengen}, is given by equation $f(x_0,x_1,x_2)=0$, where $f$ is a homogeneous polynomial of degree $d$.
Viewed as a hypersurface of $V_{\ell}$, $D$ is cut out on $V_{\ell}$ by a quadric with equation $Q(x^I)_{|I|=\ell}=0$, where $Q$ is a 
homogeneous polynomial of degree 2 in the variables $\{x^I\}_{|I|=\ell}$. 

The \emph{double plane} associated with $D$ is the double cover $\psi: \mathbb{D}^* \to\PP^2$ branched along $D$. 
It can be embedded in $\PP(1,1,1,\ell)$ as a hypersurface $\mathbb{D}^*_a$ defined by a (weighted homogeneous) equation of the form $az^2=f(x_0,x_1,x_2)$, for any $a\in \mathbb C^ *$.  Under the identification of $\PP(1,1,1,\ell)$ with $W_{\ell}$, we see that $\mathbb{D}^ *_a$ is cut out on $W_{\ell}$ by a quadric in $\PP^ {N_\ell+1}$ of the form $az^ 2=Q(x^I)_{|I|=\ell}$. 

Consider the sublinear system $\mathcal Q_{\ell}$ of $|\cO_{W_{\ell}}(2)|$ of surfaces cut out on 
$W_{\ell}$ by the quadrics of $\PP^ {N_{\ell}+1}$ with equation of the form $az^ 2=Q(x^I)_{|I|=\ell}$. 

When $a \neq 0$, the quadrics in question are such that their polar hyperplane with respect to $P$ has equation $z=0$. When $a=0$,  
such a quadric is singular at $P$, it represents the cone with vertex at $P$ over the quadric in $\Pi =\{z=0\} \cong \PP^ {N_{\ell}}$ with equation $Q(x^I)_{|I|=\ell}=0$ and it cuts out on $W_\ell$ a cone, with vertex at 
$P$, over a quadric section of $V_\ell$.

In particular, $\mathcal Q_{\ell}$ contains the codimension 1 sublinear system $\mathcal Q^c_{\ell}$ of all such cones with vertex at $P$ over a quadric section of $V_\ell$, thus $\dim(\mathcal Q_{\ell})=N_{d}+1$. Moreover $\mathcal Q_{\ell}$ is stable under 
the action of $\mathbb G_{\ell}$ on $W_{\ell}$. 

We  set $\widetilde {\mathcal Q}_{\ell}:=\rho^*(\mathcal Q_{\ell})$, which is a sublinear system of $|\mathcal O_{Z_\ell}(2)|$.  Note that $\widetilde {\mathcal Q}_{\ell}$ contains the sublinear system $\widetilde {\mathcal Q}^ c_{\ell}=\rho^*(\mathcal Q^c_{\ell})$ of all divisors 
of the form $2E$ plus a divisor in $\mathcal L_\ell(d)$.  

We denote by $\mathcal Q^*_{\ell}$ the dense open subset of $\mathcal Q_{\ell}$ of points corresponding to smooth surfaces. Since no surface in $\mathcal Q^*_{\ell}$ passes through $P$, we may and will identify $\mathcal Q^*_{\ell}$ with its pull--back via $\rho$ on $Z_{\ell}$, which is a dense open subset of $\widetilde {\mathcal Q}_{\ell}$ sitting in the complement of $\widetilde {\mathcal Q}^ c_{\ell}$. We denote by $\mathcal{I\!\!D}^*_{\ell}\to \mathcal Q^*_{\ell}$ the universal family.  

A surface  $\mathbb{D}^ *\in \mathcal Q^*_{\ell}$ cuts out on $V_{\ell}$ a smooth curve $D\in \mathcal U_{\ell}$ and 
conversely; indeed the projection from $P$ realizes $\mathbb{D}^*$  as the double cover of $\PP^ 2$ branched along $D$.  This yields  a  surjective morphism 
\[
\mathcal Q^*_{\ell}\ni \mathbb{D}^ *\stackrel{\beta}{\longrightarrow} \mathbb{D}^ *\cap V_{\ell} := D \in\mathcal U_{d}\,,
\] which sends the double plane $\mathbb{D}^*$ to its branching divisor $D$.  
This morphism is equivariant under the actions of $\mathbb G_{\ell}$ on both $\mathcal Q^*_{\ell}$ and $\mathcal U_{d}$, where  $\mathbb G_{\ell}$ acts on $\mathcal U_{d}$ via the natural action of the quotient group ${\rm PGL}(3,\mathbb C)$, see \eqref {eq:gr}. The morphism $\beta$ is not injective, its fibers being  isomorphic to $\mathbb C^ *$.

As an immediate consequence of  Lemma \ref{lem:canclass}, we  have:
\blem\label{lem:candegree} Let   $D$ be a smooth curve in $\PP^2$ of even degree $d=2\ell \geqslant 2$, let $\psi: \mathbb{D}^ * \to\PP^2$ be the double cover branched along $D$, let
$\Gamma\subset \PP^2$ be a projective curve of degree $m$ not containing  $D$, and let $\Gamma^ *$ be its pull--back via $\psi$ considered as a curve in $Z_{\ell}$. One has
\be\la{eq:intersection-index} K_{\ell}\cdot \Gamma^*=-m(d+6)\,.\ee
\elem

In the setting of  Lemma \ref {lem:candegree}, consider the diagram
\be\label{eq:diagr}
\xymatrix{
  &&C^*\ar[d]_{\psi'}  \ar[rr] ^ {\nu^ *} &&\,\, \Gamma^ *\ar[d]^\psi \\
 &&  \,\, C  \ar[rr]_{\nu}&& \Gamma && 
 }
\ee
where $\nu$ and $\nu^ *$  are the normalization morphisms and $\psi$ and $\psi'$ have degree 2 (to ease notation, here 
we have identified $\psi$ with its restriction to $\Gamma^*$). 

Let $\delta:=\delta(D,\Gamma)$. It could be that $C^ *$ splits into two components both isomorphic to $C$; in this case $\delta=0$. If $\delta=0$ and the genus of $C$ is zero, then $C^*$ certainly splits.
 Suppose that $C^*$ is irreducible, and
let $g$ and $g^*$ be the geometric genera of $\Gamma$ and $\Gamma^ *$ (i.e. the arithmetic genera of $C$ and $C^*$, respectively). Since $\psi'$ has exactly $\delta$ ramification points, the Riemann-Hurwitz formula yields
\be\label{eq:g-star} 2(g^*-1)=4(g-1)+\delta\,.\ee


\subsection{The odd degree case}\label{ss:odd} Fix  a  line $h \in |\cO_{\PP^2}(1)|$, and 
let  $D$ be a smooth curve in $\PP^2$ of degree $d = 2\ell -1 \geqslant 1$, which intersects $h$ transversely. 
We denote by $\mathcal U^ h_{d}$ the open subset of $\mathcal U_{d}$ consisting of such curves. 

For each $D\in \mathcal U^ h_{d}$, we consider the reducible curve of degree $d+1=2\ell$
$$\Delta := D + h \in |\cO_{\PP^2}(d+1)|$$  and 
the double cover $\psi \colon \mathbb{D}^ * \to\PP^2$, branched along $\Delta$. The difference with the even degree case is that $\mathbb{D}^*$ is no longer smooth, but it has double points at the $d$ points in $D\cap h$. 
In any event, as in the even degree case, we can consider the set $\mathcal Q^ *_{\ell; h}\subset |\mathcal O_{Z_{\ell}}(2)|$ of all such surfaces $\mathbb{D}^ *$, with its universal family $\mathcal{I\!\!D}^*_{\ell; h}\to \mathcal Q^{*}_{\ell; h}$ which parametrizes all double planes $\psi\colon \mathbb{D}^ * \to\PP^2$ as above. We still have the morphism
$$\beta: \mathcal Q^{*}_{\ell; h}\to \mathcal U^ h_{d}$$associating to $\mathbb{D}^*$ the branching divisor $\Delta$ of $\varphi\colon \mathbb{D}^ * \to\PP^2$ minus $h$. 

The group acting here is no longer the full group $\mathbb G_{\ell}$ but its subgroup  $\mathbb G_{\ell;h}$ which fits in the exact sequence
$$1\to \mathbb C^ *\to \mathbb G_{\ell;h}\to {\rm Aff}(2,\mathbb C)\to 1,$$
where ${\rm Aff}(2,\mathbb C)$ is the \emph{affine group} of all projective transformations in  ${\rm PGL}(3,\mathbb C)$ stabilizing $h$.

Keeping the setting and notation of \S \ref {ss:even}, Lemma  \ref {lem:candegree} still holds, as well as diagram \eqref {eq:diagr}. If $\Gamma$ intersects $h$ at $m$ distinct points which are off $D$, then the double cover $\psi': C^ *\to C$ has $\delta+m\geqslant m>0$ ramification points, where $\delta = \delta(D,\Gamma)$ as above. In particular, $C^*$ is  irreducible, and \eqref {eq:g-star} is replaced by 
\be\label{eq:odd}
2(g^*-1)=4(g-1)+\delta+m.
\ee


\section {The main result}\label{sec:result}

In this section we prove the following:

\begin{theorem}\label{prop:even} Let $\delta \geqslant 0$ be an integer such that, for a very general curve $D$ in $\PP^2$ of  degree $d=2\ell-\epsilon$, where $\epsilon\in \{0,1\}$, there exists an integral curve $\Gamma \subset \PP^2$, $\Gamma \neq D$, of geometric genus $g$  and degree $m$ 
with  $\delta(D,\Gamma)=\delta$. Then
\begin{equation}\label{eq:principal-ineq}
4g+\delta\geqslant m(d+2\epsilon-8)+5.
\end{equation}
\end{theorem}

The proof of Theorem \ref{prop:even} will be done in \S \ref {ss:proof}. First we need some more preliminaries, which we collect in the next subsection.  We keep all notation and conventions introduced so far.

\subsection{Constructing appropriate families}\la{sit:stategy} Fix  integers $m \geqslant 1$ and $g \geqslant 0$. Let $\mathcal{H}$ be the locally closed subset of $\mathcal L_m$, whose
points correspond to integral curves $\Gamma \subset \PP^2$ of degree $m$ and  geometric genus $g$; $\mathcal{H}$ is a quasiprojective variety. 
We let $\mathcal U \to \mathcal{H}$ be the universal curve.
 
\subsubsection{The even degree case}\label{ssec:even} Fix an even positive integer $d=2\ell$ and a non-negative integer $\delta$.  Consider the locally closed subset  $\mathcal I$ of  $\mathcal{H}\times  \mathcal Q^*_{\ell}$  of pairs 
$(\Gamma,\mathbb{D}^ *)$ such that $\Gamma$ does not coincide with the branch curve $D$ of $\psi: \mathbb{D}^ *\to \PP^ 2$ and 
$\delta(D,\Gamma)=\delta$. Remember that we may equivalently interpret $\mathbb{D}^ *$ as a surface in $W_\ell$ or in $Z_\ell$. Each irreducible component of $\mathcal I$ is fixed by the obvious action of  $\mathbb G_{\ell}$ on  $\mathcal{H}\times \mathcal Q^*_{\ell}$.

For any $(\Gamma,\mathbb{D}^ *)\in \mathcal I$, the pull--back $\Gamma^* \subset \mathbb{D}^ *$ of $\Gamma$ via $\psi$ is a reduced curve in $Z_{\ell}$.  Hence there is a morphism 
$\mu: \mathcal I\to \mathcal K$, where $\mathcal K$ is the Hilbert scheme of curves of $Z_{\ell}$. We let $\mathcal V\to \mathcal K$ be the corresponding universal family.  The map $\mu$ is equivariant under the actions of $\mathbb G_{\ell}$ on both $\mathcal I$ and $\mathcal K$.

Let $\pi_1\colon \mathcal I\to \mathcal{H}$ and $\pi_2\colon  \mathcal I\to  \mathcal Q^*_{\ell}$ be the two projections. 
Under the hypotheses of Theorem \ref {prop:even} and with notation as in \S\;\ref{ss:even}, the following holds.

\begin{lem}\label{ass:1} There exists an irreducible component $I$ of $\mathcal I$  which dominates $\mathcal Q_{\ell}$ via $\pi_2$. Hence $I$ dominates also $\mathcal U_{d}\subset \mathcal L_{d}$ via $\beta\circ \pi_2$. 
\end{lem}

Given $I$ as in Lemma \ref {ass:1}, we choose an irreducible, smooth subvariety $\mathcal B$ of $I$, such that $\pi_2$ restricts to an \'etale morphism of $\mathcal B$ onto its image, which is dense in $\mathcal Q_{\ell}$. To place our objects in the context of \S \ref{App}, consider  the universal family $\mathcal{I\!\!D}^ *_{\ell}\to \mathcal Q^{*}_{\ell}$ (cf.\;\S\;\ref{ss:even})  of double planes $\mathbb{D}^*$ [resp.\ $\mathcal V\to \mathcal K$ of curves $\Gamma^*\subset \mathbb{D}^*$]. Up to possibly shrinking $\mathcal B$ and performing an \'etale cover of it, the morphisms $\mathcal B \stackrel{\pi_2}{\longrightarrow} \mathcal Q^*_{\ell}$ and $\mathcal B \stackrel{\mu}{\longrightarrow} \mathcal K$ give rise to families
\[ 
\mathcal D:= \pi_{2}^ *(\mathcal{I\!\!D}^*_{\ell})\stackrel{p}\longrightarrow  \mathcal B \quad\mbox{and}\quad \mathfrak C:=\mu^ *(\mathcal V)\stackrel{\mathfrak q} \longrightarrow  \mathcal B\,.\]over $\mathcal B$ fitting in diagram \eqref {eq:in1}. We may assume  that there exists a simultaneous normalization  $\nu$ and a family $\mathcal C \stackrel{q}{\longrightarrow} \mathcal B$ as in  
\eqref {eq:tess}, with $\mathcal C$ smooth fitting in \eqref{eq:bidg},  where $X=Z_{\ell}$.

\subsubsection{The odd degree case}\la{ss:odc} Fix now an odd positive integer $d=2\ell-1$ and a non-negative integer $\delta$, and fix a line $h$ in $\PP^ 2$. We consider the locally closed subset $\mathcal I$ 
of $\mathcal{H}\times  \mathcal Q_{\ell;h}$ consisting  of pairs 
$(\Gamma,\mathbb{D}^ *)\in \,\mathcal{H}\times \mathcal Q^*_{\ell}$ such that
$\Gamma$ is not contained in the branch divisor $\Delta$ of $\psi: \mathbb{D}^ *\to \PP^ 2$, 
$\delta(D,\Gamma)=\delta$, and $\Gamma$ intersects $h$ at $m$ distinct points which are off $D$. 

For any point $(\Gamma,\mathbb{D}^ *)\in \mathcal I$, the pull--back $\Gamma^* \subset \mathbb{D}^ *$ of $\Gamma$ via $\psi$ is an integral curve in $Z_{\ell}$.  So, we still have the morphisms 
$\mu: \mathcal I\to \mathcal K$,  $\pi_1\colon \mathcal I\to \mathcal{H}$ and $\pi_2\colon \mathcal I\to  \mathcal Q^ *_{\ell;h}$ equivariant under 
actions of $\mathbb G_{\ell;h}$. 

As before, we have the following

\begin{lem}\label{ass:2} There exists an irreducible component $I$ of $\mathcal I$  which dominates $\mathcal Q^ *_{\ell;h}$ via $\pi_2$. 
\end{lem}

As in the even case, given $I$ as in Lemma \ref {ass:2}, we may construct a smooth $\mathcal B$ having an \'etale, dominant morphism to $\mathcal Q_{\ell;h}$, together with families

\[ 
\mathcal D:= \pi_{2}^ *(\mathcal{I\!\!D}^*_{\ell,h})\stackrel{p}\longrightarrow  \mathcal B, \quad \mathfrak C:=\mu^ *(\mathcal V)\stackrel{\mathfrak q} \longrightarrow  \mathcal B,\,\]
fitting in  diagram \eqref {eq:in1}. Consider a simultaneous normalization $\nu$ and a family $\mathcal C \stackrel{q}{\longrightarrow} \mathcal B$ as in  
\eqref {eq:tess}, with $\mathcal C$ smooth. In view of Lemmata \ref {ass:1} and \ref {ass:2}, the constructed families fit in diagram \eqref {eq:bidg}, with $X=Z_{\ell}$.

\bigskip

In both cases, the next lemma allows to apply Proposition \ref{prop:CL-3}  in our setting.

 \blem\label{lem:upper-bound} For any $d>0$, 
 $\mathcal C\stackrel{q} \longrightarrow\mathcal B$ is a covering family, i.e.  $z({\mathcal{C}}) =3$.
 \elem
\bproof By the discussion in \S \ref {sit:stategy}, for $d$ even $\phi(\mathcal C)$ is stable under the action of $\mathbb G_{\ell}$ on $Z_\ell$, which is transitive; for $d$ odd $\phi(\mathcal C)$ is stable under the action of $\mathbb G_{\ell;h}$, which is transitive on the dense open subset of $Z_{\ell}$ whose complement is $\pi^ {-1}(h)\cup E$. This proves the assertion.
\eproof

\subsection{Proof of Theorem \ref{prop:even}}\label {ss:proof} 

Our proof follows the one of Theorem (1.2) in \cite {CL}. First we recall the following useful fact.

\blem[See {\cite [Lemma (3.1)]{CL}}] \la{lem:us} Let $g: V\to W$ be a linear map of finite dimensional vector spaces. Suppose that $\dim(g(V))>k$. Let $\{V_i\}_{i\in I}$ be a family of vector subspaces of $V$, such that $\bigcup_{i\in I}V_i$ spans $V$, and
for any pair $(i,j)\in I\times I$, with $i\neq j$, there is a finite sequence $i_1=i, i_2,\ldots, i_{t-1},i_t=j$ of distinct elements of $I$ with
$\dim(g(V_{i_h}\cap V_{i_{h+1}}))\geqslant k$, for all $h\in \{1,\ldots, t-1\}$. Then there is an index $i\in I$ such that $\dim(g(V_i))>k$.
\elem

\subsubsection{The even degree case} \label {ss:proofeven} 
We need to construct a suitable  subfamily of $\mathfrak C\to \mathcal B$ with the  covering property. 

Fix  a general point $b_0\in \mathcal B$, and let $\Gamma_0^ *$ and $\mathbb{D}_0^ *$ be the corresponding elements of the families $\mathfrak C\to \mathcal B$ and 
$\mathcal D\to \mathcal B$, respectively. 

Let $\mathcal L$ be the open subset of the linear system $\mathcal L_{\ell}(d-1)$ as in \eqref{eq:lel} consisting of the surfaces $F\in \mathcal L_{\ell}(d-1)$ which do not contain $\Gamma^ *_0$. A general such surface $F$ meets $\Gamma_0^ *$ transversely. By genericity, we may suppose that all surfaces $F$ defined by the pull--back via $\pi: Z_{\ell} \to V_{\ell}\cong \PP^2$ of degree $d-1$ monomials in the variables $x_0,x_1, x_2$ belong to $\mathcal L$.
For a given $F\in \mathcal L$, we denote by $\mathcal B_F$ the subvariety of $\mathcal B$ parameterizing  all double planes in  $\mathcal D\to \mathcal B$ containing the complete intersection curve of $F$ and $\mathbb{D}_0^ *$. In addition, for a general point $\xi\in \Gamma_0^ *$ we let $\mathcal B_{F,\xi}$ denote the subvariety of $\mathcal B_F$ parameterizing  all surfaces in  $\mathcal D\to \mathcal B_F$ which pass through $\xi$.

\blem\label{lem:ea} For $F\in \mathcal L$ and $\xi\in \Gamma_0^ *$  as above one has
\[
\dim(\mathcal B_F)= 3\,\,\, \text{and}\,\,\, \dim (\mathcal B_{F,\xi})= 2\,.
\]
Furthermore, $b_0$ is a smooth point of both $\mathcal B_F$ and $\mathcal B_{F,\xi}$.
\elem

\begin{proof} Consider the sublinear system $\Lambda_F$ of $\widetilde {\mathcal Q}_{\ell}=\rho^*( {\mathcal Q}_{\ell})$ on $Z_\ell$ consisting of all surfaces containing the complete intersection curve $F\cap \mathbb{D}_0^*$. Imposing to the surfaces in $\Lambda_F$ the condition to contain a general point of $F$, the divisor $F+2E$ splits off, and the residual surface sits in $\mathcal L_\ell(1)$. Hence $\Lambda_F$ contains a  codimension 1 sublinear system consisting of surfaces of the form $2E+F+L$, with $L$ varying in $\mathcal L_{\ell}(1)$, which has dimension 2. Hence $\dim (\Lambda_F)=3$. Since $\mathcal B$ dominates $\mathcal Q_{\ell}$ via $\pi_2$, which is finite on $\mathcal B$, and $\mathcal B_F$ is the inverse image of $\Lambda_F$, one has $\dim(\mathcal B_F)=\dim (\Lambda_F)= 3$.
The proof is similar for $\mathcal B_{F,\xi}$. The final assertion follows by the genericity assumptions.
\end{proof} 

We denote by $T_0$ the tangent space to $\mathcal B$ at $b_0$, and by $T_F$ and $T_{F,\xi}$ the 3 and 2--dimensional subspaces of $T_0$ tangent to $\mathcal B_F$ and to $\mathcal B_{F,\xi}$ at $b_0$, respectively. 

\blem\label{lem:gene} One has:
 \begin{enumerate}
\renewcommand\labelenumi{\rm (\alph{enumi})}
\renewcommand\theenumi{\rm (\alph{enumi})}
\item
$\bigcup_{F\in \mathcal L} T_F$ spans $T_0$;
 \item 
given $F\in \mathcal L$, the union $\bigcup_{\xi\in \Gamma_0^*} T_{F,\xi}$ spans $T_F$.
 \end{enumerate}
\elem

\begin{proof} (a) Since $\pi_2$ is \'etale on $\mathcal B$, $T_0$ is isomorphic to the tangent space to $\mathcal Q_{\ell}$ at $\mathbb{D}_0^ *$. Remember that, by \S \ref {ss:even},  the double plane $\mathbb{D}_0^ *$, considered in $W_{\ell}$, is cut out by a quadric with equation $z^ 2=Q(x^I)_{|I|=\ell}$. So $T_0$
can be identified with the vector space of homogeneous quadratic polynomials of the form $az^ 2-G(x^I)_{|I|=\ell}$ modulo the one-dimensional linear space spanned by $z^ 2-Q(x^I)_{|I|=\ell}$ and by the linear space of quadratic polynomials in $\{x^I\}_{|I|=\ell}$ defining $V_{\ell}$. \footnote{An explanation is in order. Consider a vector space $V$ and a nonzero vector $v\in V$, along with the associated projective space $\PP(V)$ and the corresponding point $[v]\in\PP(V)$. Then the tangent space $T_{[v]}\PP(V)$ can be canonically identified with ${\rm Hom}(\langle v\rangle, V/\langle v \rangle)\cong V/\langle v \rangle$.} Hence $T_0$ can be identified with the vector space of quadratic polynomials in $\{x^I\}_{|I|=\ell}$,  modulo the vector space of quadratic polynomials in $\{x^I\}_{|I|=\ell}$ defining $V_{\ell}$. This, in turn,  can be identified  with the vector space $S_{d}$ of homogeneous polynomials of degree $d$ in $x_0,x_1,x_2$.

Now $T_F$ can be identified with the  vector subspace  $S_{d}(f)$ of $T_0\cong S_{d}$ of polynomials of the form $fh$, where $f$ is a fixed homogeneous polynomial of degree $d-1$ (determined by $F$), and $h$ is any linear form. By assumption on $\mathcal L$, 
$\bigcup_{F\in \mathcal L} T_F$ contains all monomials of degree $d$, which do span $S_{d}$. 

(b)  Given $F$,  $T_{F,\xi}$ can be identified with the vector space of homogeneous polynomials of the form $fh$,  where $h$  vanishes at $\pi(\xi)\in \mathbb P^ 2$. These polynomials  do span $T_F\cong S_{d}(f)$. 
 \end{proof}

Next we consider the restrictions 
\[
\mathcal D_F\stackrel{p}\longrightarrow  \mathcal B_F, \quad \mathfrak C_F \stackrel{\mathfrak q} \longrightarrow  \mathcal B_F,\,\,\, \text{and}\,\,\, \mathcal D_{F,\xi}\stackrel{p}\longrightarrow  \mathcal B_{F,\xi}, \quad \mathfrak C_{F,\xi} \stackrel{\mathfrak q} \longrightarrow  \mathcal B_{F,\xi},\]respectively, of the families \[\mathcal D\stackrel{p}\longrightarrow  \mathcal B\,\,\,\text{ and}\,\,\, \mathfrak C\stackrel{\mathfrak q} \longrightarrow  \mathcal B.\]

\bprop\la{prop:cov} For general $F\in \mathcal L$ and  $\xi\in \Gamma_0^ *$,  the  families 
\[
\mathfrak C_F \stackrel{\mathfrak q} \longrightarrow  \mathcal B_F\,\,\, \text{and}\,\,\, \mathfrak C_{F,\xi} \stackrel{\mathfrak q} \longrightarrow  \mathcal B_{F,\xi}
\]
have the covering property.
\eprop

\begin{proof}  We prove the assertion for $\mathfrak C_F \stackrel{\mathfrak q} \longrightarrow  \mathcal B_F$. The proof for $\mathfrak C_{F,\xi} \stackrel{\mathfrak q} \longrightarrow  \mathcal B_{F,\xi}$ is similar (and analogous to the proof of the  corresponding  statement in \cite [Theorem (1.2)] {CL}), hence it can be left to the reader. 

Let $\mathbb M$ be the set of all monomials of degree $d-1$ in $x_0,x_1,x_2$.  Consider  the family $\{F_M\}_{M\in \mathbb M}$, where $F_M\in \mathcal L$ is defined by the pull--back via $\pi\colon Z_\ell\to V_\ell\cong\PP^2$ of the monomial $M$. Take two monomials $M',M''$ which differ only in degree 1, i.e., their lowest common multiple $U$ has degree $d$. Then $\mathcal B_{F_ {M'}}\cap  \mathcal  B_{F_{M''}}$ contains the pull--back via $\pi_2$ of an open, dense subset of the pencil $\langle \mathbb{D}_t^* \rangle$ spanned by $\mathbb{D}_0^*$ and $F_U$, where $F_U$ is the pull--back via $\pi$ of the monomial $U$. The base locus of this pencil does not contain $\Gamma_0^ *$. Therefore, $\Gamma_0^ *$ varies in a non-trivial one-parameter family  $\langle \Gamma_t^ *\rangle$ together with members $\mathbb{D}_t^*$ varying in the pencil $\langle \mathbb{D}_t^* \rangle$. 

Next we apply Lemma \ref {lem:us} with 
\begin{itemize}
\item $V=T_0$;
\item $W=H^ 0(\Gamma_0^ *,N_{\Gamma_0^ *|Z_\ell})$;
\item the linear map $g$ induced by the characteristic map (see Proposition \ref {prop:CL} (a));
\item  the family of subspaces $\{V_i\}_{i\in I}$ given by 
$\{T_{F_M}\}_{M\in \mathbb M}$. 
\end{itemize}
For each pair of monomials $M',M''$, there is a sequence of monomials $M_1=M',M_2,\ldots, M_{t-1}, M_t=M''$, such that for all $i=1,\ldots, t-1$, the lowest common multiple of $M_i$ and $M_{i+1}$ has degree $d$. The above argument implies that  $g({T_{F_{M_i}}} \cap {T_{F_{M_{i+1}}}} )$  has dimension at least 1, for all $i=1,\ldots, m-1$. Furthermore, one has $\dim (g(T_0))\geqslant 2$, because $\mathfrak C\to \mathcal B$ is a covering family (see (b) of Proposition \ref {prop:CL} and Lemma \ref{lem:upper-bound}). By Lemma \ref {lem:us} there is a monomial $M \in \mathbb M$ such that $\dim(g(T_{F_M}))\geqslant 2$; by virtue of Lemma \ref{lem:ea}, this implies that $\mathfrak C_{F_M}\to \mathfrak B_{F_M}$ is a covering family. This proves the assertion.\end{proof}

To finish the proof of  Theorem \ref {prop:even} in this case, consider the covering family $\mathfrak C_{F,\xi} \stackrel{\mathfrak q} \longrightarrow  \mathcal B_{F,\xi}$, with $F\in \mathcal L$ and $\xi\in \Gamma_0^ *$ general. Using \eqref{eq:foc-set-deg}, \eqref{eq:intersection-index}, and
\eqref {eq:g-star}, for $b=(\Gamma^*_b,\mathbb{D}_b^*)\in \mathcal B_{F,\xi}$ general (see \S \ref{ssec:even}) we deduce
\begin{equation}\la{eq:ff}
\deg (\Phi_b^{\rm sm}) \leqslant 4(g-1)+\delta+2m(\ell+3)= 4(g-1)+\delta+m(d+6)\,.
\end{equation}
On the other hand, by construction and by (a) of Proposition \ref {prop:CL-3}, 
\begin{equation}\la{eq:fff}
\deg (\Phi_b^{\rm sm}) \geqslant 1+\deg (\Gamma_b^*\cap F)=1+2(d-1)m\,.
\end{equation}
Comparing \eqref{eq:ff} and \eqref{eq:fff} gives \eqref{eq:principal-ineq}.

\subsubsection{The  odd degree case}  \label {ss:proofodd} The proof runs exactly as in the case of even $d$, so we will be brief and leave the details to the reader. 

Fix  again $b_0\in \mathcal B$, $\Gamma_0^ *$ and $\mathbb{D}_0^ *$ as in the even degree case. Following what we did in \S \ref{ss:odc}, we replace $D_b$ by $D_b+h$, where $h\subset\PP^2$ is  a general line. In the present setting we let $\mathcal L$ be the open subset of $\mathcal L_{\ell}(d)$ consisting of the surfaces $F\in\mathcal L_{\ell}(d)$ which do  not contain $\Gamma^*_0$. Again we may assume that all surfaces $F$ defined by the pull--back via $\pi$ of degree $d$ monomials in the variables $x_0,x_1, x_2$ belong to $\mathcal L$. Given $F\in \mathcal L$, we define $\mathcal B_F$ and $\mathcal B_{F,\xi}$ as in the even degree case, and the analogue of Lemma \ref {lem:ea} still holds.  Then, with the usual meaning for  $T_0$,  $T_F$ and $T_{F,\xi}$, the analogue of Lemma \ref {lem:gene} holds.  Similarly as in Proposition \ref {prop:cov}, the covering property holds for the restricted families

\[
\mathcal D_F\stackrel{p}\longrightarrow  \mathcal B_F, \quad \mathfrak C_F \stackrel{\mathfrak q} \longrightarrow  \mathcal B_F,\quad \text{and}\quad \mathcal D_{F,\xi}\stackrel{p}\longrightarrow  \mathcal B_{F,\xi}, \quad \mathfrak C_{F,\xi} \stackrel{\mathfrak q} \longrightarrow  \mathcal B_{F,\xi}\,.\]
We conclude finally as in the even degree case: \eqref {eq:ff} holds with no change, whereas \eqref{eq:fff} has to be replaced by
$$
\deg (\Phi_b^{\rm sm}) \geqslant 1+\deg (\Gamma_b\cap F)=1+2dm\,,
$$and  again, \eqref{eq:principal-ineq}  follows.

\section {Genera of subvarieties: a survey}\label{sec:appendix}
As we mentioned in the Introduction, inequality \eqref{eq:us} 
allows to bound the genera of curves in double planes from below. In this section we provide  a brief survey on genera bounds for subvarieties in different type of varieties, and discuss several conjectures. All varieties are supposed to be projective, reduced, irreducible, and defined over $\CC$. The geometric genus $p_g(Y)$ of a variety $Y$ is the geometric genus of a smooth model of $Y$. 

Two important sources of interest in bounding genera are: the Clemens Conjecture on count of rational curves in Calabi-Yau varieties (\cite{Cl2}), and the celebrated Kobayashi Conjecture on hyperbolicity of general hypersurfaces in $\PP^n$ of sufficiently large degree (\cite{Kob}). Recall (\cite{DR}, \cite{Kob}) that the Kobayashi hyperbolicity of a variety $X$ implies the algebraic hyperbolicity, and in particular, absence of rational and elliptic curves in $X$. A part concerning the Clemens Conjecture started with
the following theorem.

\begin{theorem}\label{Clemens} {\rm (H. Clemens \cite{Cl1})} The geometric genera of curves in a  very general hypersurface  $X$ of degree $d \geqslant 2n-1$ in $\PP^n$ satisfy the inequality $g\geqslant 
\frac{1}{2} (d - 2n + 1) + 1$.
\end{theorem}

This shows, in particular, the absence of rational curves in  very general surfaces of degree $d\geqslant 5$ in $\PP^3$. One of the subsequent results in higher dimensions  was 

\begin{theorem} {\rm (E.~Ballico \cite{bal})}  There is an effective bound $\varphi(n)$ such that a very  general  hypersurface of degree $d\geqslant \varphi(n)$ 
in $\PP^n$ is algebraically hyperbolic. 
\end{theorem}

A better effective bound was provided by Geng Xu (\cite{Xu1}).
For instance (\cite{CZ}), it follows from the results in \cite{Xu1} that a general sextic threefold in $\PP^4$ is algebraically hyperbolic. 

The  Demailly algebraic hyperbolicity theorem states the following.

 \begin{theorem}\label{thm:Demailly} {\rm (J.-P.~Demailly \cite{Dem1})}
For any hyperbolic subvariety $X\subset\PP^n$ there exists $\epsilon>0$ such that, for any curve $C\subset X$,  the geometric genus $g$ of $C$ is bounded below in terms of the degree: $g\geqslant \epsilon \deg(C)+1$. Therefore, the curves of bounded genera in $X$ form bounded families. 
\end{theorem}

Due to a recent proof of the Kobayashi Conjecture, Theorem \ref{thm:Demailly} can be applied to general (in Zariski sense) hypersurfaces in $\PP^n$.

\begin{theorem} {\rm (D.~Brotbek  \cite{Bro},  Y.-T.~Siu \cite{Siu})}
A general hypersurface of sufficiently large degree in $\PP^n$ is Kobayashi hyperbolic. 
\end{theorem}

For effective estimates of degrees of hyperbolic hypersurfaces, see Y.~Deng (\cite{Deng,Deng1}); see also J.-P.~Demailly \cite{Dem3} for a survey and a simplified proof. 

L.~Ein obtained some analogs of Clemens' estimate in higher dimensions. 
 
\begin{theorem} {\rm (L.~Ein \cite{Ein1, Ein2})} Let $M$ be a smooth projective variety of dimension $n\geqslant 3$, let $L\to M$ be
an ample line bundle, and let $X \in |dL|$ be a very general member. Then for $d \geqslant 2n - \ell$ any subvariety $Y\subset X$ of dimension $\ell$ has  positive geometric genus, and for 
$d \geqslant 2n - \ell + 1$,  $Y$ is of general type.
\end{theorem}

In the case $M=\PP^n$ there is a sharper bound.

\begin{theorem} {\rm (C.~Voisin \cite{Voi1,Voi2})} Let $X$ be a  very general hypersurface of degree $d \geqslant 2n- \ell - 1$ in $\PP^n$, $n\geqslant 3$. Then  for 
$d \geqslant 2n -�� l + 1$ any subvariety $Y\subset X$ of dimension $l\le n-3$ has  positive geometric genus, and for $d \geqslant 2n -�� l$, $Y$ is of general type.
\end{theorem}

Sharper bounds are known also for certain toric varieties (A.~Ikeda \cite{Ik}). 
For $M=\PP^n$ a further improvement is due to G.~Pacienza. 

\begin{theorem} {\rm (G.~Pacienza \cite{Pac2})} For $n \geqslant 6$ and for a  very general hypersurface $X\subset\PP^n$ of degree $d \geqslant 2n-2$, any subvariety $Y\subset X$ is of general type.
\end{theorem}

Geng Xu improved Ein's theorem as follows.

\begin{theorem} {\rm (G.~Xu \cite{Xu5})} Let $X$ be a  very general complete intersection of $m \leqslant n - 3$ hypersurfaces of degrees $d_1,\ldots,d_m$ in $\PP^n$, where $d_i \geqslant 2$ $\forall i$, and let $Y\subset X$ be a reduced and irreducible divisor. Let $d=d_1 +\cdots+d_m$. Then
$p_g(Y)\geqslant n-1$ if $d \geqslant n+2$, and
$Y$ is of general type if $d > n + 2$.
\end{theorem}

See also Geng Xu (\cite{Xu4}), C.~Chang and Z.~Ran (\cite{ChR}), L.~Chiantini, A.-F.~Lopez, and Z.~Ran (\cite{CLR}), H.~Clemens and Z.~Ran (\cite{CR}),  S.S.-T.~Lu and Y.~Miyaoka (\cite{LM1}), and L.-C.~Wang (\cite{LCWang1, LCWang2}). The results in \cite{LCWang1} include some classes of divisors in Calabi-Yau hypersurfaces of degree $d=n+1$ in $\PP^n$.

Let us mention several sporadic results. See also, e.g., R.~Beheshti (\cite{Beh}), M.~Bernardara (\cite{Ber}), L.~Bonavero and A.~Hoering (\cite{BH}), T.D.~Browning and P.~Vishe (\cite{BV}), I.~Coskun and E.~Reidl (\cite{CoRe}), O.~Debarre (\cite{De}), K.~Furukawa (\cite{Fu1, Fu2}),  J.~Harris, M.~Roth, and J.~Starr (\cite{HRS}),  J.~Kollar (\cite{Kol}). Concerning the Clemens Conjecture on rational curves in quintic threefolds and Mirrow Symmetry, see, e.g.,  M.~Kontsevich (\cite{Kon}), A.~Libgober and J.~Teitelbaum (\cite{LT}), D.A.~Cox and S.~Katz (\cite{CK}), T.~Coates and A.~Givental (\cite{CoGi})  and the references therein.  

\begin{theorem} \begin{itemize}
\item {\rm (G.~Pacienza \cite{Pac1}, E.~Riedl and D.~Yang \cite{RY})} Let $X \subset \PP^n$
be a very general  hypersurface of degree $d$. If either $n=6$ and $d=2n-3$, or $n\geqslant 7$ and $\frac{3n+1}{2} \le d \le
2n - 3$, then $X$ contains lines but no other rational curves.
\item {\rm (D.~Shin \cite{Shin2})} A general hypersurface of degree $d > \frac{3}{2} n - 1$ in $\PP^n$ does not contain any smooth conic; however {\rm (S.~Katz \cite{Ka})}, a general quintic threefold in $\PP^4$ does.  
\item {\rm (S.~Katz \cite{Ka}, P.~Nijsse \cite{Ni}, T.~Johnsen and S..L.~Kleiman \cite{JK1},  J.~D'Almeida \cite{DA}, E.~Cotterill \cite{Co0, Co1}, E.~Ballico and C.~Fontanari \cite{BaFo}; cf.\ also E.~Ballico \cite{bal1} and A.~L.~Knutsen \cite{Kn1, Kn2})} A general
quintic threefold $X$ in $\PP^4$ contains only finitely many rational curves of degree $\le 12$, and each rational curve 
$C$ of degree $\le 11$ either is smooth and embedded in $X$ with a balanced normal bundle $\mathcal{O}(-1) \oplus \mathcal{O}(-1)$, or is a plane $6$-nodal quintic.  
\item {\rm (D.~Shin \cite{Shin1}, G.~Mostad~Hana and T.~Johnsen \cite{HJ}, E.~Cotterill \cite{Co2})} A general hypersurface  of degree $7$ in $\PP^5$ does not contain any rational curve of degree $d\in\{2,\ldots,16\}$. 
\item {\rm (G.~Ferrarese and D.~Romagnoli \cite{FR})}  The degree of an elliptic curve on a very general hypersurface $X$ of degree $7$ in $\PP^4$ is a multiple of $7$. 
\item {\rm (B.~Wang \cite{Wang2})} A general hypersurface of degree $54$
in $\PP^{30}$ does not contain any rational quartic curve.
\item {\rm (B.~Wang \cite{Wang1})} A  very general  hypersurface of
degree $d \geqslant 2n - 1$ in $\PP^n$ does not contain any smooth elliptic curve.
\end{itemize}
\end{theorem}
For $n=3$, Geng Xu replaced Clemens' initial genus bound in Theorem \ref{Clemens} by the  optimal one. See also L.~Chiantini and A.F.~Lopez (\cite{CL}) for an alternative proof and some generalizations.

\begin{theorem}\label{Xu} {\rm (G.~Xu \cite{Xu3})} The genera of curves on a very general surface of degree $d \geqslant 5$ in $\PP^3$ satisfy the inequality $g \geqslant \frac{1}{2}d(d - 3) - 2$, and this bound is sharp.
For $d \geqslant 6$ this sharp bound can be achieved only by a tritangent hyperplane section.  
\end{theorem}

Let  ${\rm Gaps}(d)$
be the set
 of all non-negative integers which cannot be realized as geometric genera of irreducible curves on a very general surface of degree 
$d$ in $\PP^3$. This set is union of finitely many disjoint and separated integer intervals. 
By Xu's Theorem \ref{Xu}, the first gap interval is
${\rm Gaps}_0(d) = [0, d(d-3)/
2 - 3]$. For $d=5$, this is the only gap interval. For $d\geqslant 6$,
the next gap interval is ${\rm Gaps}_1(d)= [
d(d-3)/2+
2, d^2 - 2d- 9]$ (\cite{CFZ1}). One can show (\cite{CFZ2}) that $\max ({\rm Gaps}(d))=O(d^{8/3})$. The latter is based on certain existence results. For arbitrary smooth (not necessarily general) surfaces in $\PP^3$, we have the following existence result.

\begin{theorem}\label{smooth-surfaces} {\rm (\cite{JAChen}, \cite{CFZ2})}  
There exists a function $c(d)\sim d^3$ such that, for any smooth surface $S$ in  $\PP^3$ of degree $d$
and any $g \geqslant c(d)$, $S$ carries
a reduced, irreducible nodal curve of geometric genus $g$, whose nodes can be prescribed generically on $S$. 
\end{theorem}

To formulate an analog in higher dimensions, we recall the following notion.
Let $Y$ be an irreducible variety of dimension $s$. A singular point $y \in Y$ is called an \emph{ordinary singularity of multiplicity} $m$ ($m > 1$), if the Zariski tangent space of $Y$ at $y$ has dimension $s + 1$, and
the (affine) tangent cone to $Y$ at $y$ is a cone with vertex $y$ over a smooth hypersurface of degree $m$ in $\PP^s$.

The next result was first established by J.A.~Chen (\cite{JAChen}) for curves in $n$-dimensional varieties, and then in \cite{CFZ2} for subvarieties of arbitrary dimension $s\leqslant n-1$ \footnote{We are grateful to J.A.~Chen for pointing out his nice paper \cite{JAChen} that we ignored when writing \cite{CFZ2}. We apologize for our ignorance.}.

\begin{theorem}\label{higher-dim} {\rm (\cite{JAChen}, \cite{CFZ2})} Let $X$ be an  irreducible, smooth, projective variety of dimension $n>1$, let $L$ be a very ample divisor on $X$, and let $s\in\{1,\ldots,n-1\}$. 
Then there is an  integer $p_{X,L,s}$ 
such that 
for any $p\geqslant p_{X,L,s}$ one can find an irreducible complete intersection $Y=D_1\cap\ldots\cap D_{n-s}\subset X$ of dimension $s$ with at most ordinary points of multiplicity $s+1$ as singularities such that $p_g(Y)=p$, where $D_i\in |L|$ for $i=1,\ldots,n-s-1$ are smooth and transversal and $D_{n-s}\in |mL|$ for some $m\geqslant1$. Moreover, for $n\geqslant 3$ and $s=1$ one can find a \emph{smooth} curve $Y$ in $X$ of a given genus $g(Y)=p\geqslant p_{X,L,1}$. 
 \end{theorem}
 
Recall the famous:

\smallskip

{\bf Green-Griffiths-Lang Conjecture.}  (\cite{GG, La}; see also \cite{Berc}, \cite{Dem2}) \emph{Let $X$ be a projective variety of general type. Then there exists a proper closed subset $Z\subset X$ such that any subvariety $Y\subset X$ not of general type is contained in $Z$. }

\smallskip

The following conjecture is inspired by the previous results and by the Green-Griffiths-Lang Conjecture in the surface case. 
 
{\bf Conjecture.} \emph{There exists a strictly growing function $\varphi(d)$, with natural values, such that the set of curves of geometric genus $g\leqslant\varphi(d)$ in any smooth surface $S$ of degree $d\geqslant 5$ in $\PP^3$ is finite. }
 
 \smallskip
 
Notice (\cite{CFZ2}) that for any $g\geqslant 0$ and $d\geqslant 1$ one can find a smooth surface $S\subset \PP^3$ of degree $d$ carrying a nodal curve of genus $g$. 
 Notice also that any smooth quartic surface in $\PP^3$ contains an infinite countable set of rational curves, hence the restriction $d\geqslant 5$ is necessary. Let us mention a few facts supporting the conjecture. According to B.~Segre (\cite{Seg}) the number of lines on a smooth surface of degree $d\geqslant 3$ does not exceed $(d-2)(11d-6)$. The celebrated Bogomolov theorem (\cite{Bo}) says that the number of rational and elliptic curves on a surface of general type with $c_1^2>c_2$ is finite. Moreover,  due to Y.~Miyaoka, this number admits a uniform estimate:

 \begin{theorem} {\rm (Y.~Miyaoka \cite{Mi1})}
Let $S$ be a minimal smooth projective surface of general type satisfying the inequality for Chern numbers $c_1^2>c_2$.
Then the number of irreducible curves of genus $0$ and $1$ on $S$ is bounded by a function
of $c_1$ and $c_2$. 
 \end{theorem}
 
 Analogous facts are true under certain weaker assumptions on Chern numbers (Y.~Miyaoka \cite{Mi2}), or on the singularities of rational and elliptic curves in $S$ (S.S.-Y.~Lu and Y.~Miyaoka \cite{LM2}). It is plausible that the number of curves of genus $g\leqslant \varphi(d)$ on a smooth surface of degree $d$ in $\PP^3$ can be uniformly bounded above by a function of $d$. The  conjecture above is coherent with the following ones.
 
 {\bf Conjecture.} {\rm (C.~Voisin \cite{Voi3})} \emph{ Let $X\subset\PP^n$ be a very general hypersurface of degree $d\geqslant n+2$. Then the degrees of rational curves in $X$ are bounded.}
 
 \smallskip
 
 {\bf  Conjecture.}  {\rm (P.~Autissier, A.~Chambert-Loir, and C.~Gasbarri \cite{ACLG})} \emph{ Let $X$ be a smooth projective variety of general type with  the canonical line bundle $K_X$. Then there exist real numbers $A$ and $B$, and a proper Zariski closed subset $Z \subset X$ such that for any curve $C$ of geometric genus $g$ in $X$ not contained in $ Z$, one has
$\deg_C(K_X) \leqslant A(2g - 2) + B$.}
 


\end{document}